\documentclass[11pt]{article}
\usepackage{amssymb}
\usepackage{amsmath}
\usepackage{theorem}
\usepackage{latexsym}
\reversemarginpar
\theoremheaderfont{\bf} %\theoremstyle{break}
\theorembodyfont{\sl}
\newtheorem{theo}{Theorem}[section]
{\theorembodyfont{\rm} \newtheorem{defi}[theo]{Definition}}
{\theorembodyfont{\rm} \newtheorem{exa}[theo]{Example}}
{\theorembodyfont{\rm} \newtheorem{rem}[theo]{Remark}}
\newtheorem{prop}[theo]{Proposition}

\newtheorem{lemma}[theo]{Lemma}
{\theorembodyfont{\rm}}
{\theorembodyfont{\rm}}
\newenvironment{proof}{{\sc Proof:}}{\mbox{}\hfill$\Box$\par}
\newcommand{\eqnref}[1]{~\mbox{$(${\rm \ref{#1}}$)$}}

\newcommand{\junk}[1]{}
\newcommand{\DS}{\displaystyle}

\newcommand{\N}{{\mathbb N}}
\newcommand{\F}{{\mathbb F}}

\newcommand{\cC}{{\mathcal C}}

\newcommand{\cS}{{\mathcal S}}

\newcommand{\ve}[1]{\mbox{$\varepsilon^{(#1)}$}}

\newcommand{\rank}{\mbox{\rm rank}\,}
\newcommand{\ord}{\mbox{\rm ord}}
\newcommand{\spann}{\mbox{\rm span}\,}

\newcommand{\Char}{\mbox{\rm char}}
\newcommand{\AutF}{\mbox{${\rm Aut}_{\mathbb F}$}}
\newcommand{\im}{\mbox{\rm im}\,}

\newcommand{\dist}{\mbox{\rm dist}}
\newcommand{\wt}{\mbox{\rm wt}}

\newcommand{\floor}[1]{\mbox{$\lfloor{#1}\rfloor$}}

\newcommand{\ideal}[1]{\mbox{$\langle{#1}\rangle$}}
\newcommand{\z }{\pmb{z}}
\newcommand{\Azs}{\mbox{$A[\z;\sigma]$}}

\newcommand{\p}{\mbox{$\mathfrak{p}$}}
\renewcommand{\v}{\mbox{$\mathfrak{v}$}}

\newcommand{\lideal}[1]{\mbox{$^{^{\bullet\!\!}}\langle{\, #1\, }\rangle$}}

\newcounter{abc}
\newcounter{def}
\newenvironment{romanlist}{\begin{list}{(\roman{abc})\hfill}{\usecounter{abc}
     \topsep-.5ex \labelwidth.7cm \leftmargin.7cm \labelsep0cm
     \rightmargin0cm \parsep0ex \itemsep.6ex
     \partopsep1ex}}{\end{list}}
\newenvironment{alphalist}{\begin{list}{(\alph{abc})\hfill}{\usecounter{abc}
     \topsep-.5ex \labelwidth.7cm \leftmargin.7cm \labelsep0cm
     \rightmargin0cm \parsep0ex \itemsep.6ex
     \partopsep1ex}}{\end{list}}
\newenvironment{arabiclist}{\begin{list}{(\arabic{abc})\hfill}{\usecounter{abc}
     \topsep-.5ex \labelwidth.7cm \leftmargin.7cm \labelsep0cm
     \rightmargin0cm \parsep0ex \itemsep.6ex
     \partopsep1ex}}{\end{list}}

\title{On doubly-cyclic convolutional codes}
\author{Heide Gluesing-Luerssen\footnote{
       University of Groningen, Department of Mathematics, P.~O.~Box 800,
       9700 AV Groningen, The Netherlands; gluesing@math.rug.nl}
       \ and Wiland Schmale\footnote{
       Department of Mathematics, University of Oldenburg,
       26\,111 Oldenburg, Germany;
       schmale@mathe\-matik.uni-oldenburg.de}
       }
\date{October 13, 2004}

\begin{document}
 \maketitle

\begin{abstract}
\noindent
Cyclicity of a convolutional code (CC) is relying
on a nontrivial automorphism of the algebra
$\F[x]/(x^n-1)$, where $\F$ is a finite field. If this
automorphism itself has certain specific cyclicity
properties one is lead to the class of doubly-cyclic
CC's.
Within this large class Reed-Solomon and BCH convolutional codes can be defined.
After constructing doubly-cyclic CC's, basic properties are derived
on the basis of which distance properties of Reed-Solomon convolutional codes
are investigated.
This shows that some of them are optimal or near optimal with
respect to distance and performance.
\end{abstract}



%%%%%%%%%%%%%%%%%%%%%%%%%%%%%%%%%%%
\section{Introduction}
\setcounter{equation}{0}
%%%%%%%%%%%%%%%%%%%%%%%%%%%%%%%%%%%

Despite the fact that convolutional codes are as important for applications
as block codes, their mathematical description is much less
developed, and there has been growing activity to fill this gap during the last
decade, see, e.~g.,
\cite{RoSm99,SGR01,PPS04,HRS03,GS03,GS04,GL03}.

The gap in the mathematical theory of block and convolutional codes is particularly
big when it comes to the notion of cyclicity.
Cyclic convolutional codes (shortly, cyclic CC's or just CCC's)
have been introduced and investigated by Piret and Roos in~\cite{Pi76,Ro79};
for definitions see below.
Their approach has much later been extended in~\cite{GS04} to a
theoretical framework which exhibits many features in close analogy to the
well known theory of cyclic linear block codes.
It turned out that the class of CCC's contains plenty of
codes with very good performance and distance properties,
see also \cite{GS03,GL03}.

In this article we construct a specific subclass of CCC's
where the generating polynomial has an additional cyclic structure,
see Section~3.
Among these are Reed-Solomon type doubly-cyclic CC's, for which distance
properties are derived in Section~4. More general results, leading to BCH convolutional
codes, are indicated in Section~5.
A minimum of prerequisites can be found in Section~2.

One standard way of defining CC's is as follows.
%%%%%%%%%%%%%%%%%%%%%%%%%%%%
\begin{defi}\label{D-CC}
Let $\F$ be any finite field.
A {\em convolutional code\/} $\cC\subseteq\F[\z]^n$ with (algebraic)
parameters $(n,k,\delta)$ is a submodule of the form $\cC=\im G$, where
$G\in\F[\z]^{k\times n}$ is a right-invertible matrix such that
$\delta=\max\{\deg\gamma\mid \gamma\text{ is a }k\text{-minor of }G\}$.
We call~$G$ a {\em generator matrix\/} of the code.
The number~$n$ is called the {\em length},~$k$ is the {\em dimension}, and~$\delta$ is
called the {\em overall constraint length\/} of the code.
\end{defi}
%%%%%%%%%%%%%%%%%%%%%%%%%%%%%

By elementary matrix and module theory over $\F[\z]$ one realizes that
a CC with parameters $(n,k,\delta)$ is just a direct summand of $\F[\z]^n$ of rank~$k$
and that the overall contraint length~$\delta$ does not depend on the choice of the
generating matrix~$G$ for~$\cC$. Details can be found for instance in
\cite{Fo70,McE98,GS04}.
In the coding literature a right invertible matrix is often called
{\em basic\/}~\cite[p.~730]{Fo70} or
{\em delay-free and non-catastrophic}, see~\cite[p.~1102]{McE98}.

It is well-known that each submodule of
$\F[\z]^n$ has a minimal generator matrix in the sense of the next
definition~\cite[Thm.~5]{Fo70} or~\cite[p.~495]{Fo75}.
In the same paper~\cite[Sec.~4]{Fo75} it has been shown how to derive such a matrix
from a given generator matrix in a constructive way.
For a row vector $v\in\F[\z]^n$ we will denote by $\deg v$ the maximum degree of its
components. The zero vector has degree~$-\infty$.
%%%%%%%%%%%%%%%%%%%%%%%%%%%%%
\begin{defi}\label{D-minBasis}
Let $G\in\F[\z]^{k\times n}$ be a matrix with rank~$k$ and overall constraint
length~$\delta$ and let $\nu_1,\ldots,\nu_k$ be the degrees of the rows of~$G$.
We say that $G$ is {\sl minimal\/} if $\delta=\sum_{i=1}^k\nu_i$.
In this case the row degrees of~$G$ are uniquely determined by the submodule
$\cS:=\im G$. They are called the {\sl Forney indices\/} of~$\cS$.
The largest Forney index is called the {\sl memory\/} of~$\cS$.
\end{defi}
%%%%%%%%%%%%%%%%%%%%%%%%%%%%%
The notion ``minimal'' stems from the (simple) fact that for an arbitrary generator
matrix~$G$ one has $\delta\leq\sum_{i=1}^k \nu_i$.
Thus, in a minimal generator matrix the rows degrees have been reduced to their minimal
values.
Using such a generator matrix it is easily seen that a code with overall constraint
length zero can be regarded as a block code.

An important quality characteristic of a code is its so-called free
distance. It measures the error-correcting capability.
For a polynomial vector $v=\sum_{j=0}^N v_j \z^j\in\F[\z]^n$, where $v_j\in\F^n$,
the {\em weight\/} is defined as $\wt(v)=\sum_{j=0}^N\wt(v_j)$ where the weight
of~$w_j\in\F^n$ denotes the usual Hamming weight.
Then the {\em (free) distance\/} of a code $\cC\subseteq\F[\z]^n$ is, just like for
block codes, defined as
$\dist(\cC):=\min\{\wt(v)\mid v\in\cC,\;v\not=0\}$.

%%%%%%%%%%%%%%%%%%%%%%%%%%%%%%%%%%%%%%%%%%%%%%%%%
\section{Preliminaries for cyclic convolutional codes}
\setcounter{equation}{0}
%%%%%%%%%%%%%%%%%%%%%%%%%%%%%%%%%%%%%%%%%%%%%%%%%
As usual a cyclic block code of length $n$ and dimension $k$ over the field~$\F$
will be described as a principal ideal in the algebra $A:=\F[x]/\ideal{x^n-1}$.
We always assume that $\Char(\F)$ does not divide $n$.
We have the natural isomorphisms
\[
   \p:\F^n\rightarrow A,\
   (v_0,\ldots,v_{n-1})\mapsto\sum_{i=0}^{n-1}v_ix^i\text{ and } \v:=\p^{-1}.
\]
The weight function on~$A$ is defined such that~$\p$ is an isometry between~$A$
and $\F^n$ endowed with the usual Hamming metric,
i.~e., $\wt(a):=\wt\big(\v(a)\big)$ for all $a\in A$. Let
\begin{equation}\label{equ-pi}
x^n-1=\prod^{r-1}_{i=0}\pi_i
\end{equation}
be the prime factorization over $\F[x]$.
Since we assume $\Char(\F)$ and $n$ to be coprime,
the normed prime polynomials $\pi_i$ are all different.
According to this factorization the algebra $A$
decomposes into a direct sum of minimal cyclic block codes,
which can be generated by the (primitive) idempotents
$\ve{i}, 0\le i\le r-1$. We have
\begin{equation}\label{e-epspi}
    \ve{i}\,\text{mod}\,\pi_j=\delta_{ij}\text{ for all }i,j=0,\ldots,r-1,
\end{equation}
and their existence is guaranteed by the Chinese Remainder Theorem.
The idempotents are uniquely determined by $A$
and\eqnref{e-epspi} implies
\begin{equation}\label{equ-eps1}
\ve{i}=\beta\prod_{j\ne i}\pi_j  \text{  for some unit } \beta \in \F.
\end{equation}
The cyclic code $\ideal{\ve{i}}$ is minimal and also isomorphic to $\F[x]/\ideal{\pi_i}$
and in addition one has
\begin{equation}\label{equ-eps2}
 \dim_{\F} \ideal{\ve{i}}=\deg\pi_i.
\end{equation}
Moreover, any cyclic block code of length~$n$ over~$\F$ is generated by a sum of
idempotents, which is unique up to ordering of the summands.

In the convolutional setting,
the vector space~$\F^n$ has to be replaced by
$\F[\z]^n:=\{\sum_{\nu=0}^N\z^{\nu}v_{\nu}\mid N\in\N_0,\,v_{\nu}\in\F^n\}$ and,
consequently, the ring~$A$ by the polynomial ring
$$A[\z]:=\big\{\sum_{j=0}^N\z^ja_j\,\big|\, N\in\N_0,\,a_j\in A\big\}$$
over~$A$.
The natural extensions of the maps~$\p$ and $\v$  are given by
\begin{equation}\label{e-p}
\p\big(\sum_{\nu=0}^N\z^{\nu}v_{\nu}\big)=\sum_{\nu=0}^N\z^{\nu}\p(v_{\nu})
\text{ and } \v:=\p^{-1}
\end{equation}
where, of course, $v_{\nu}\in\F^n$ and thus $\p(v_{\nu})\in A$ for all~$\nu$.
This map is an isomorphism of $\F[\z]$-modules.
Note that $\p$ and $\v$ are isometries if we define $\wt(g):=\wt\big(\v(g)\big)$ for all
$g\in A[\z]$.

It is now tempting to define a cyclic convolutional code (CCC) to be an ideal
$A[\z]$ or more precisely, to declare a code $\cC\subseteq\F[\z]^n$ as cyclic if
$\p(\cC)$ is an ideal in $A[\z]$.
It has been shown in \cite[Thm.~3.12]{Pi76} and \cite[Thm.~6]{Ro79} that this does not
result in any codes other than block codes, see also \cite[Prop.~2.7]{GS04}.
Led by this negative result, a more general notion of cyclicity has been
introduced for convolutional codes~\cite{Pi76,Ro79,GS04}.
It makes use of an automorphism of the $\F$-algebra~$A$.
Thus, let $\AutF(A)$ to be the group of all $\F$-automorphisms on~$A$.
Detailed information on this group can be found in \cite[Sec.~3]{GS04}.
In particular, it is shown that in general there are quite a lot of
automorphisms and how to determine them.
For later use we only wish to mention that firstly, each automorphism~$\sigma\in\AutF(A)$
is uniquely determined by the value of $\sigma(x)$, and secondly,
for each~$a\in A$ such that $\ord(a)\mid n$ and each~$k\in\{0,\ldots,n-1\}$
the assignment $\sigma(x)=a^kx$ determines an automorphism.
We will mainly make use of this type of automorphism though in general there may
be many others, too.

Picking an arbitrary automorphism $\sigma\in\AutF(A)$,
a new multiplication in the $\F[\z]$-module $A[\z]$ is defined via
\begin{equation}\label{e-az}
        a\z=\z\sigma(a) \text{ for all }a\in A
\end{equation}
along with associativity and distributivity.
This turns $A[\z]$ into a non-commutative $\F[\z]$-algebra which will be denoted by $\Azs$.
We call $\Azs$ the {\it Piret algebra} (over $A$ and with respect to the
automorphism~$\sigma$).
Note that it coincides with the commutative ring $A[\z]$
if~$\sigma$ is the identity. In all other cases it is a non-commutative ring.
In particular, it is important to distinguish between left and right coefficients of~$\z$.
The coefficients can be moved to either side by applying the
rule\eqnref{e-az} since~$\sigma$ is invertible.
Multiplication inside~$A$ remains the same as before.
Hence~$A$ is a commutative subring of $\Azs$.
Due to this very specific non-commutativity the ring $\Azs$ is also called a
{\em skew-polynomial ring}.
Since $\sigma|_{\F}=\text{id}_{\F}$, the ordinary commutative polynomial ring
$\F[\z]$ is a subring of $\Azs$, too.
As a consequence, $\Azs$ inherits the (left and right) $\F[\z]$-module structure from
$A[\z]$ .
For us, only the left module structure will be important.
In particular, the map~$\p$ from\eqnref{e-p} is an isomorphism between the left
$\F[\z]$-modules $\F[\z]^n$ and $\Azs$ (notice that in~$\p$ the coefficients are on
the right of~$\z$).

Now we declare a submodule $\cC\subseteq\F[\z]^n$ to be $\sigma$-{\em cyclic\/} if
$\p(\cC)$ is a left ideal in $\Azs$.
Cyclic CC's have been investigated in detail in the papers \cite{Pi76,Ro79,GS04,GS03}
and it turned out that there are many good codes that are not block codes.
See~\cite{GS04} for more details.
In the same paper an algebraic theory of CCC's has been developed
where  in the context of Piret algebras notions like non-catastrophicity, dimension of a
code, and overall constraint length could be handled successfully.
In the next section multiple use of these results will be  made.

%%%%%%%%%%%%%%%%%%%%%%%%%%%%%%%%%%%
\section{Construction of doubly-cyclic codes}
\setcounter{equation}{0}
%%%%%%%%%%%%%%%%%%%%%%%%%%%%%%%%%%%
In this section we will give a construction of convolutional codes with
parameters $(n,k,km)$ where~$m$ is the memory.
It is based on cyclic block codes as discussed in the previous section.
The distances of a subclass of these codes will be computed in Section~4.

Let us fix an automorphism~$\sigma$.
It is easy to see that~$\sigma$ induces a permutation on the set
\[
    E=\{\ve{0},\ldots,\ve{r-1}\}.
\]
Remember that according to\eqnref{e-epspi} the $i$th idempotent corresponds to the $i$th
prime factor of $x^n-1$.
Since $\sigma(\ve{i})=\ve{j}$ implies $\deg\pi_i=\deg\pi_j$,
an automorphism can induce a nontrivial permutation on $E$ only if
the degrees of the prime factors of $x^n-1$ are not all pairwise
different.
In this case there exists a subset $S\subset E$ such that $S\cap\sigma(S)=\varnothing$.
Let from now on $\sigma\in\AutF(A)$ be such an automorphism.

We then fix a subset~$S$ and define $b\in\N$ such that
\begin{equation}\label{e-S}
    S\cap\sigma^j(S)=\varnothing\text{ for all }1\leq j\leq b.
\end{equation}
Note that this implies $\sigma^i(S)\cap\sigma^j(S)=\varnothing$ for all $0\leq i<j$
such that $j-i\leq b$.
Let $s:=|S|$. Then $(b+1)s\leq r$ and $
(b+1)s=r\Longleftrightarrow E=\bigcup_{i=0}^{\,b}\sigma^i(S)$.

Consider now the cyclic block code generated by
\begin{equation}\label{e-c}
   c:=\sum_{\ve{i}\in S}\ve{i}.
\end{equation}
It is the direct sum of the minimal block codes $\ideal{\ve{i}}$
and based on\eqnref{equ-eps1} and\eqnref{equ-eps2} one obtains
\begin{equation}\label{e-k}
  k:=\dim_{\F}\ideal{c}=\sum_{\ve{i}\in S}\deg\pi_i.
\end{equation}
A basis is, for instance, given by the elements $c,\,xc,\ldots,x^{k-1}c$.
Equation\eqnref{e-S} can now also be expressed via the orthogonality
\begin{equation}\label{e-corthog}
   \sigma^i(c)\sigma^j(c)=0\text{ for all }0\leq i<j\text{ such that }j-i\leq b.
\end{equation}
%%%%%%%%%%%%%%%%%%%%%%%%%%%%%%%%%%%%%%%%%%
\begin{exa}\label{E-SandG}
\begin{alphalist}
\item  Let $q=4$, $n=15$, and $\alpha$ be a primitive element for~$\F$.
       We compute
        \begin{align*}
           x^{15}-1=&\,(x+1)(x+\alpha^2)(x+\alpha)(x^2+\alpha^2x+1)(x^2+\alpha
                      x+\alpha)(x^2+x+\alpha^2)\\
                    &\,(x^2+\alpha^2x+\alpha^2)(x^2+\alpha x+1)(x^2+x+\alpha)
        \end{align*}
        and order the idempotents $\ve{0},\ldots,\ve{8}$ according to the ordering of
    the factors.
        For the automorphism we consider~$\sigma$ defined by $\sigma(x)=\alpha x$.
        Then one can show that the permutation $\sigma|_{E}:E\longrightarrow E$ has
    the cycles
    \[
      \big(\ve{0},\ve{1},\ve{2}\big)\big(\ve{3},\ve{4},\ve{5}\big)
      \big(\ve{6},\ve{7},\ve{8}\big).
    \]
        Let now, for instance, $S=\{\ve{0},\ve{3},\ve{6}\}$ then
        $S\cup \sigma(S)\cup\sigma^2(S)$ is a disjoint union and is equal to
        $E=\{\ve{0},\dots ,\ve{8}\}$, the full set of idempotents.
\item  Let $q=2$ and $n=31$.  One computes
        \begin{align*}
            x^n-1=&\,(x+1)(x^5+x^4+x^3+x^2+1)( x^5+x^2+1)(x^5+x^4+x^3+x+1)\\
                  &\,(x^5+x^3+x^2+x+1)( x^5+x^3+1)(x^5+x^4+x^2+x+1)\ .
        \end{align*}
        Let the idempotents $\ve{0},\dots ,\ve{6}$ be numbered accordingly.
        In this situation the assignment $\sigma(x):=x^3$ leads to an automorphism
        with $\sigma(\ve{0})=\ve{0}$ and $\sigma(\ve{k})=\ve{k+1}$ for $1\le k\le 5$.
        Now, for instance, defining $S=\{\ve{1},\ve{4}\}$ one has $|S|=2$ and
        \[
        \{\ve{1},\dots ,\ve{6}\}=S\cup \sigma(S)\cup\sigma^2(S)
        \]
        as a disjoint union. This example is typical in some sense, since $x-1$ is always one
        of the prime factors of $x^n-1$.
        Thus, if $x^n-1$ has no further linear factors, then $S$ can
        only contain idempotents different from $ \ve{0}$.

\end{alphalist}

\end{exa}
%%%%%%%%%%%%%%%%%%%%%%%%%%%%%%%%%%%%%%%%%%%%%%
The following example introduces CCC's of Reed-Solomon type which will
be further investigated in later sections.
%%%%%%%%%%%%%%%%%%%%%%%%%%%%%%%
\begin{exa}\label{E-RScodes}
Let $\F=\F_q$ be a field of size~$q$ and let $n:=q-1$.
Furthermore let $\alpha\in\F$ be a primitive element, thus $\ord(\alpha)=n$.
Then the prime factor decomposition of $x^n-1$ is given by
$x^n-1=\prod_{i=0}^{n-1}\pi_i$, where $\pi_i=x-\alpha^i$.
We pick $k\in\N$ such that $1\leq k\leq\frac{n}{2}$ and choose $\sigma\in\AutF(A)$
such that
\[
    \sigma(x)=\alpha^k x\ .
\]
Since $\ord(\alpha^k)\mid n$ this does indeed define an automorphism on~$A$.
Since $\ve{j}=\beta_j\prod_{i\ne j}(x-\alpha_i)$ for $0\le j\le n-1$ and some
$\beta_j\in\F^{\ast}$, the automorphism~$\sigma$ acts on the idempotents as follows:
\[
   \sigma(\ve{j})=\ve{j-k\;\text{mod}\;n}, \text{ for } j=0,\dots ,n-1.
\]
Define
\[
  S:=\{\ve{n-k},\dots ,\ve{n-1}\} \text{ and }
  b:=\floor{\frac{n}{k}}-1.
\]
Then Equation\eqnref{e-S} is satisfied
 and  due to the restriction $k\leq\frac{n}{2}$ we have $b\geq 1$.
Let now~$c$ be as in\eqnref{e-c}.
Then $\ideal{c}$ is a $k$-dimensional cyclic block code with generator polynomial
\[
  f:=\prod_{l=0}^{n-k-1}(x-\alpha^{\,l})\in\F[x]
\]
and~$k$ is as in\eqnref{e-k}.
This shows that $\ideal{c}$ is a Reed-Solomon code of length~$n$.
It is well-known, see e.~g. \cite[Thm.~6.6.2]{vLi99},
that
\[
       \dist\ideal{c}= n-k+1.
\]
\end{exa}
%%%%%%%%%%%%%%%%%%%%%%%%%%%%%%%%%%%%%%%%%%%

We return now to the general situation and introduce what will be called a
{\em doubly-cyclic convolutional code}.
Using the ingredients from\eqnref{e-S} --\eqnref{e-k} along with the
automorphism~$\sigma$ and the isomorphism from\eqnref{e-p} we define the
matrix
\begin{equation}\label{e-G}
  G:=\sum_{\nu=0}^m \z^{\nu}G_{\nu}\in\F[\z]^{k\times n} \text{ where }
  G_{\nu}:=\begin{pmatrix}
              \v\big(\sigma^{\nu}(c)\big)\\
              \v\big(\sigma^{\nu}(xc)\big)\\
                \vdots\\
              \v\big(\sigma^\nu(x^{k-1}c)\big)
            \end{pmatrix}\in\F^{k\times n},
\end{equation}
and where
\begin{equation}\label{e-m}
    1\leq m\leq b.
\end{equation}

The matrix~$G$ above might look artificial.
However, it becomes quite natural once considered over the appropriate Piret algebra
$\Azs$. Recall that $\F[\z]^n\cong\Azs$ as left $\F[\z]$-modules via the isomorphism $\p$
in\eqnref{e-p} and also recall the skew multiplication defined via\eqnref{e-az}.
Define
\begin{equation}\label{e-g}
   g:=c\sum_{\nu=0}^m\z^{\nu}=\sum_{\nu=0}^m\z^{\nu}\sigma^{\nu}(c)\in\Azs.
\end{equation}
Then we obtain $x^ig=\sum_{\nu=0}^m\z^{\nu}\sigma^{\nu}(x^ic)$ for all $i\in\N_0$
and, due to left $\F[\z]$-linearity of~$\v$,
\[
  G=\begin{pmatrix}\v(g)\\\v(xg)\\ \vdots\\ \v(x^{k-1}g)\end{pmatrix}.
\]
In Theorem~\ref{T-alg} we will show that $\im G=\v(\lideal{g})$ where
$\lideal{g}:=\{fg\mid f\in\Azs\}$ is the left ideal generated by~$g$.
Moreover, we will see that~$G$ is right invertible and thus defines a cyclic
convolutional code. Also dimension and overall constraint length of this code are derived.
For a subclass of doubly-cyclic CC's we will compute distances and extended row distances
in the next section.
%%%%%%%%%%%%%%%%%%%%%%%%%%%%%%%%%%%%%%%%%%%%%%%%%%%%%%%%%%%%%%%%%%%%
\begin{theo}\label{T-alg}
Let the data be as in\eqnref{e-S} --\eqnref{e-g}.
Then
\begin{alphalist}
\item $g:=c\big(1+\z\sigma(c)\big)\big(1+\z\sigma^2(c)\big)\cdot\ldots\cdot
               \big(1+\z\sigma^m(c)\big)$.
\item We have $gu=c$ where
      $u=\big(1-\z\sigma^m(c)\big)\big(1-\z\sigma^{m-1}(c)\big)\cdot\ldots\cdot
       \big(1-\z\sigma(c)\big)$.
      Furthermore, $u$ is a unit in $\Azs$ and
      $u=1-\z\big(\sigma(c)+\ldots+\sigma^m(c)\big)$.
\item Define $\cC:=\im G$. Then $\cC=\v\big(\lideal{g}\big)$.
      Thus, $\cC$ is a cyclic submodule of $\F[\z]^n$. Moreover, $\rank\cC=k$.
\item $\cC$ is a cyclic convolutional code, or in other words, a direct summand of
      $\F[\z]^n$.
      Equivalently, the matrix~$G$ is right invertible.
\item The matrix $G$ is minimal in the sense of Definition~\ref{D-minBasis}.
\item $\cC$ is a code with parameters $(n,k,km)$ and memory~$m$.
      In particular, all Forney indices of the code~$\cC$ are equal to~$m$.
\end{alphalist}
The convolutional code~$\cC$ will be called a doubly-cyclic code.
\end{theo}
%%%%%%%%%%%%%%%%%%%%%%%%%%%%%%%%%%%%%%%%%%%%%%%%%%%%%%%%%%%%%%%%%%%%
\begin{proof}
(a) We proceed by induction. For $m=1$ we have, since $\sigma(c)$ is idempotent,
\[
   g=c+\z\sigma(c)=c+\z (\sigma(c))^2=c+c\,\z\sigma(c)=c(1+\z\sigma(c)).
\]
Let now
$\sum_{\nu=0}^{m-1}\z^{\nu}\sigma^{\nu}(c)
  =c\big(1+\z\sigma(c)\big)\big(1+\z\sigma^2(c)\big)\cdot\ldots\cdot
    \big(1+\z\sigma^{m-1}(c)\big)$.
Then
\[
   \sum_{\nu=0}^{m-1}\z^{\nu}\sigma^{\nu}(c)\big(1+\z\sigma^{m}(c)\big)
  = \sum_{\nu=0}^{m-1}\z ^{\nu}\sigma^{\nu}(c)
     +\sum_{\nu=0}^{m-1}\z ^{\nu}\sigma^{\nu}(c)\z\sigma^{m}(c)
  =\sum_{\nu=0}^{m} \z ^{\nu}\sigma^{\nu}(c).
\]
The last identity follows from the fact that
\[
    \z^{\nu}\sigma^{\nu}(c)\,\z\sigma^{m}(c)=\z^{\nu+1}\sigma^{\nu+1}(c)\sigma^{m}(c)
    =\left\{\begin{array}{ll}
    0, &\text{ if }\nu<m-1\\ \z^m\sigma^m(c),&\text{ if }\nu=m-1
    \end{array}\right.
\]
due to $m\leq b$ and\eqnref{e-corthog}.
\\
(b) The equation $gu=c$ as well as the fact that~$u$ is a unit follow from~(a) along with
\[
    \big(1-\z\sigma^{\nu}(c)\big)\big(1+\z\sigma^{\nu}(c)\big)
    =\big(1+\z\sigma^{\nu}(c)\big)\big(1-\z\sigma^{\nu}(c)\big)=1,
\]
which in turn is a consequence of $\sigma^{\nu+1}(c)\sigma^{\nu}(c)=0$,
see\eqnref{e-corthog}.
The last part of~(b) can easily be shown as in~(a).
\\
For the assertions~(c) --~(f) we first have to show that the polynomial~$g$ is reduced
in the sense of \cite[Def.~4.9(b)]{GS04}.
We have
\[
   \ve{i}{g}=\left\{\begin{array}{ll}
              0 &\text{ if }\ve{i}\not\in S\\
              \sum_{\nu=0}^m\z^{\nu}\sigma^{\nu}(\ve{i})&\text{ if }\ve{i}\in S.
              \end{array}\right.
\]
This shows that the polynomials $\ve{i}g,\,\ve{i}\in S$, all have degree~$m$ and
their highest coefficients do not divide each other in~$A$ proving the reducedness of~$g$
in the above mentioned sense.
Now, application of \cite[Thm.~7.8]{GS04} yields~(c) while (d) follows from
\cite[Prop.~7.10]{GS04} along with part~(b) above.
\\
(e) To see minimality of~$G$, observe that the leading coefficient matrix is given by
\[
    \begin{pmatrix}\v\big(\sigma^m(c)\big)\\\v\big(\sigma^m(xc)\big)\\ \vdots\\
                   \v\big(\sigma^m(x^{k-1}c)\big)\end{pmatrix}.
\]
This matrix has full row rank since, by choice of~$c$, the polynomials
$c,\,xc,\ldots,x^{k-1}c$ are linearly independent in the $\F$-vector space~$A$.
Hence~$G$ is a minimal matrix due to \cite[p.~495]{Fo75}.
\\
(f) is a consequence of the previous results.
\end{proof}
Notice that by construction doubly-cyclic codes are always proper convolutional codes,
i.~e., codes with nonzero memory.
They are determined by the cyclic block code~$\ideal{c}$ and the cyclic behavior of
the automorphism~$\sigma$.

Note that part~(a) and the first statement of~(b) in Theorem~\ref{T-alg}
still remain true for $m=b+1$.
For the statements in~(c) to~(f) no restriction for~$m$ is necessary.
Later on, however,\eqnref{e-m} will be an essential assumption in order to obtain precise
informations on the free distance of doubly-cyclic codes.
We will also need the following information on various block
codes which appear in our construction.

%%%%%%%%%%%%%%%%%%%%%%%%%%%%%%%%%%%
\begin{prop}\label{P-coeffcode}
Let~$G$ and $G_{\nu}$ be as in\eqnref{e-G}
%%Heide: eingefuegt, der Konsistenz wegen
and\eqnref{e-m}.
%%endHeide
For $0\leq\mu\leq\nu\leq m$ define the matrix
\[
   G_{\mu,\nu}:=\begin{pmatrix}G_{\mu}\\ G_{\mu+1}\\ \vdots\\G_{\nu}\end{pmatrix}\in
   \F^{(\nu-\mu+1)k\times n}.
\]
and put $\cC_{\mu,\nu}:=\im G_{\mu,\nu}$.
Then $\cC_{\mu,\nu}$ is a cyclic block code given by
\[
   \cC_{\mu,\nu}
   =\ideal{\sigma^{\mu}(c)}+\ldots+\ideal{\sigma^{\nu}(c)}
   =\ideal{\sigma^{\mu}(c)}\oplus\ldots\oplus\ideal{\sigma^{\nu}(c)}.
\]
Moreover, $\dim\cC_{\mu,\nu}=(\nu-\mu+1)k$ and~$\cC_{\mu,\nu}$ has idempotent generator
\[
   {\displaystyle \sigma^{\mu}(c)+\ldots+\sigma^{\nu}(c)=
  \!\!\sum_{\varepsilon\in\sigma^{\mu}(S)\cup\ldots\cup\sigma^{\nu}(S)}
  \hspace*{-2em}\varepsilon}\hspace*{1.5em}.
\]
\end{prop}
%%%%%%%%%%%%%%%%%%%%%%%%%%%%%%%%%%%%%%%%%%
\begin{proof}
As for the first identity, observe that each code $\ideal{\sigma^i(c)}$ has dimension~$k$
and is generated by the elements $\sigma^i(c),\sigma^i(xc),\ldots,\sigma^i(x^{k-1}c)$.
This follows easily from the case $i=0$ and the fact that~$\sigma$ is an automorphism.
Therefore,
\begin{align*}
  \cC_{\mu,\nu}&=\spann_{\F}\{
        \sigma^{\mu}(c),\sigma^{\mu}(xc),\ldots,\sigma^{\mu}(x^{k-1}c),\ldots,
        \sigma^{\nu}(c),\sigma^{\nu}(xc),\ldots,\sigma^{\nu}(x^{k-1}c)\}\\
     &=\ideal{\sigma^{\mu}(c)}+\ldots+\ideal{\sigma^{\nu}(c)}.
\end{align*}
The second identity follows from\eqnref{e-corthog} along with the
inequalities $\mu\leq\nu\leq m\leq b$, see also \cite[Thm.~6.4.3]{vLi99}.
As a consequence we obtain $\dim\cC_{\mu,\nu}=(\nu-\mu+1)k$.
The form of the idempotent generator is a consequence of the fact that each
$\sigma^i(c)$ is the idempotent generator of the corresponding code.
Hence the direct sum is generated by the sum of these generators, see again
\cite[Thm.~6.4.3]{vLi99}.
\end{proof}

In special cases one can even obtain simple formulas for the distances of the codes
$\cC_{\mu,\nu}$.
As we will see next this is, for instance, the case in the situation of
Example~\ref{E-RScodes}.
%%%%%%%%%%%%%%%%%%%%%%%%%%%%
\begin{lemma}\label{L-coeffcodedist}
Let $\F$ and~$n$, the automorphism~$\sigma$ and the set~$S$ be as in
Example~\ref{E-RScodes}. 
%%Heide: auch hier die zweite Gleichung eingefuegt
Define the matrix~$G$ as in\eqnref{e-G},\eqnref{e-m} 
%%endHeide
and let the code
$\cC_{\mu,\nu}$ be as in Proposition~\ref{P-coeffcode}. Then
\[
   \dist(\cC_{\mu,\nu})=n-(\nu-\mu+1)k+1.
\]
for all $0\leq\mu\leq\nu\leq m$.
\end{lemma}
%%%%%%%%%%%%%%%%%%%%%%%%%%%%%%%
\begin{proof}
First notice that~$\sigma$ is an isometry, i.e., $\wt(a)=\wt(\sigma(a))$ for all $a\in A$.
Thus it suffices to show the result for $\mu=0$, see also Proposition~\ref{P-coeffcode}.
In the case under consideration we have
$\sigma^i(S)=\{\ve{n-(i+1)k},\ve{n-(i+1)k+1},\ldots,\ve{n-ik-1}\}$. Thus
\[
  S\cup\sigma(S)\cup\ldots\cup\sigma^{\nu}(S)=
  \{\ve{i}\mid i=n-(\nu+1)k,\ldots,n-1\}.
\]
Thus, Proposition~\ref{P-coeffcode} shows that
\[
  \cC_{0,\nu}=\Big\langle\sum_{i=n-(\nu+1)k}^{n-1}\ve{i}\Big\rangle
   =\Big\langle\prod_{i=0}^{n-(\nu+1)k-1}\pi_i\Big\rangle.
\]
Since $\pi_i=x-\alpha^i$, the generator polynomial has exactly $n-(\nu+1)k$ consecutive
powers of~$\alpha$ as zeros, proving that $\dist(\cC_{0,\nu})\geq n-(\nu+1)k+1$.
Using $\dim(\cC_{0,\nu})=(\nu+1)k$ from Proposition~\ref{P-coeffcode}
together with the Singleton bound completes the proof.
\end{proof}

%%%%%%%%%%%%%%%%%%%%%%%%%%%%%%%%%%
\section{Distance parameters for Reed-Solomon convolutional codes}
\setcounter{equation}{0}
%%%%%%%%%%%%%%%%%%%%%%%%%%%%%%%%%%%
In this section we will consider only the situation of Example~\ref{E-RScodes}.
We will compute the distances of the codes of this type and also derive lower bounds
for the extended row distances.

We begin with presenting the following upper bound on the distance of 
%%Heide: Tippfehler
convolutional
%%endHeide
codes with given algebraic parameters.
It will later provide us with some insight into the quality of the
codes constructed in the foregoing sections.
For one-dimensional codes we will see that our codes attain the generalized Singleton
bound \cite[Thm.~2.2]{RoSm99}
\begin{equation}\label{e-MDS1}
   \dist(\cC)\leq n(m+1) \text{ for any code }\cC\text{ with parameters }(n,1,m).
\end{equation}
For codes of bigger dimension we will compare the distance with the upper bound
given next.

%%%%%%%%%%%%%%%%%%%%%%%%%%%%%%%%%%%%%
\begin{prop}\label{P-estim}
Let $n=q-1$ and $\cC\subseteq\F[\z]^n$ be an $(n,k,km)_q$-code with memory~$m$ and
dimension $k>1$ and such that the memory satisfies $m\leq\frac{n}{k}-1$.
Then $\dist(\cC)\leq (m+1)(n-k+1)+(k-2)m$.
\end{prop}
%%%%%%%%%%%%%%%%%%%%%%%%%%%%%%%%%%%%%
\begin{proof}
This follows easily by using the Griesmer bound, see~\cite[Thm.~3.4]{GS03}.
Indeed, the case $i=1$ in the Griesmer bound shows that the distance~$d$ of~$\cC$
satisfies
$\sum_{l=0}^{k-1}\lceil\frac{d}{(n+1)^l}\rceil\leq n(m+1)$.
Suppose now that $d\geq(m+1)(n-k+1)+(k-2)m+1=(m+1)(n+1)-k-2m+1$.
Then the above implies
\begin{align*}
  &(m+1)(n+1)-k-2m+1+\frac{(m+1)(n+1)-k-2m+1}{n+1}\hspace*{4cm}\\
  &\hspace*{4cm}+\sum_{l=2}^{k-1}\Big\lceil\frac{(m+1)(n-k+1)-k-2m+1}{(n+1)^l}\Big\rceil\leq n(m+1),
\end{align*}
and, using that the upper floors in the sum are all at least~$1$, we obtain
$\frac{(m+1)(n+1)-k-2m+1}{n+1}\leq m$.
Hence $\frac{k+2m-1}{n+1}\geq 1$.
But this implies $m>\frac{n-k}{2}$, contradicting $m\leq\frac{n-k}{k}$ since $k\geq 2$.
\end{proof}

We will see below that in the $2$-dimensional case our codes attain this bound, hence
are optimal.
It is not clear to us whether the bound can actually be realized by a suitable
code for arbitrary dimension $k\geq2$ and memory $m\leq\frac{n}{k}-1$.

Let us repeat the situation of Example~\ref{E-RScodes}.
Thus
\begin{equation}\label{e-data}
  n:=q-1,\ 1\leq k\leq\frac{n}{2},\text{ and }
  \alpha\in\F:=\F_q\text{ such that }\ord(\alpha)=n.
\end{equation}
Then the prime factor decomposition of $x^n-1$ is given by
$x^n-1=\prod_{i=0}^{n-1}\pi_i$, where $\pi_i=x-\alpha^i$.
Choose $\sigma\in\AutF(A)$ such that
\begin{equation}\label{e-sigma}
    \sigma(x)=\alpha^k x.
\end{equation}
This assignment does indeed define an automorphism on~$A$.
It has been shown in Example~\ref{E-RScodes} that
$S:=\{\ve{n-k},\dots ,\ve{n-1}\}$ and $b:=\floor{\frac{n}{k}}-1$
satisfy\eqnref{e-S}. Thus let
\begin{equation}\label{e-c2}
   c:=\ve{n-k}+\ldots+\ve{n-1}.
\end{equation}
As shown in Example~\ref{E-RScodes}, $\ideal{c}$ is a Reed-Solomon block code.
Therefore, we call the code $\cC=\im G$ where~$G$ is in\eqnref{e-G} and\eqnref{e-m} in this
situation a {\em Reed-Solomon convolutional code}.

Below we will not only compute the (free) distance of the associated codes but also
the extended row distances.
They have been introduced in~\cite[p.~639]{ThJu83} and~\cite[p.~541]{JPB90} and are
most closely related to the performance of the code.\footnote{The row distances, as
defined in~\cite[p.~114]{JoZi99} do not give any further information.
They are all equal to the free distance $n(\delta+1)$.}
The $j$th extended row distance amounts to the minimum weight of all paths through the
state diagram starting at the zero state and which reach the zero state after exactly~$j$
steps for the first time.
In other words, it is the minimum weight of all atomic codewords of degree~$j-1$
(i.~e., length~$j$) in the sense of~\cite{McE98a}.
The details are also explained in~\cite[Sec.~3.10]{JoZi99}.
In our case where all row degrees of the matrix~$G$ are equal to~$m$
(see Theorem~\ref{T-alg}(f)), the atomic codewords are easily described.
We will confine ourselves to the following property. It follows easily from the fact
that the last~$m$ coefficient vectors of the message $u\in\F[\z]^k$ make up the
current state in the state diagram.

%%%%%%%%%%%%%%%%%%%%%%%%%%%%
\begin{rem}\label{R-molecular}
Let $G\in\F[\z]^{k\times n}$ be a minimal right-invertible generator matrix with all row
degrees equal to~$m$ and let $u\in\F[\z]^k$. Then the following are equivalent.
\begin{romanlist}
\item The codeword $uG$ is atomic (i.~e., the associated path through the state diagram
      does not pass through the zero state except for its starting and end point).
\item The polynomial $u\in\F[\z]^k$ does not have~$m$ consecutive zero coefficients
      in~$\F^k$.
\end{romanlist}
\end{rem}
%%%%%%%%%%%%%%%%%%%%%%%%%%%%%
Having this in mind, the $j$th extended row distance of the
code~$\cC=\im G$ is given by
\[
  \hat{d}^r_j:=\min\Big\{\wt(uG)\,\Big|\,
        \begin{array}{l} u\in\F[\z]^k,\,u_0\not=0,\,\deg u=j-m-1, \text{ and}\\
             \text{no~$m$ consecutive coefficients of~$u$ are zero}\end{array}\Big\}
  \text{ for all }j\geq m+1.
\]
Notice that $\deg(u)=j-m-1$ implies $\deg(uG)=j-1$ and thus the associated path has
length~$j$.
The shortest length occurring is, of course, $m+1$.
It should also be observed that in our case the extended row distances do not depend on
the choice of the minimal generator matrix~$G$.
This follows easily from the fact that, since all Forney indices are equal
to~$m$, two minimal generator matrices are related via left multiplication by some
constant regular matrix.\footnote{If not all Forney indices are identical, then in general
the extended row distances do indeed depend on the choice of the minimal generator matrix.}

Now we can formulate the result about the distance and the extended row distances of
the cyclic code under consideration.
%%%%%%%%%%%%%%%%%%%%%%%%%%%%%%%%%
\begin{theo}\label{T-dRS}
Let the data be as in\eqnref{e-data} --\eqnref{e-c2}.
Let $\cC=\im G\subseteq\F[\z]^n$ be the code with generator matrix~$G$ defined
in\eqnref{e-G} and\eqnref{e-m} where $b=\big\lfloor\frac{n}{k}\big\rfloor-1$. Then
\begin{arabiclist}
\item $\dist(\cC)=(m+1)(n-k+1)$.
\item $\hat{d}^r_j\geq (m+1)(n-k+1)+(j-1-m)(n-k(m+1)+1)$ for all $j\geq m+1$.
\end{arabiclist}
In other words, the extended row distances are bounded from below by a linear function
with slope $n-k(m+1)+1$.
\end{theo}
%%%%%%%%%%%%%%%%%%%%%%%%%%%%%%%%%
Notice that in the case $m=0$ the first part reduces to the classical result for
$k$-dimensional Reed-Solomon block codes.
Moreover, we see that for $k=1$ the codes thus constructed attain the generalized Singleton
bound\eqnref{e-MDS1}, thus are MDS codes in the sense of~\cite[Def.~2.5]{RoSm99} and
that for~$k=2$ the codes are optimal among all codes
over the same field and with the same parameters, according to Proposition~\ref{P-estim}.
For bigger~$k$ the distance stays linearly below the upper bound given in
Proposition~\ref{P-estim}.
Part~(2) shows in particular that all codewords of weight~$(m+1)(n-k+1)$ are
associated with constant messages, i.~e., messages of length~$1$.
It is worth mentioning that the slope $n-k(m+1)+1$ for the extended row distances is
optimal.
Indeed, as we will see below in\eqnref{e-vbig} for large degree the ``middle coefficients''
of a codeword are contained in the block code generated by $G_{0,m}$.
In our case this matrix has full row rank (thus no cancellation $uG_{0,m}=0$ can arise) and
the code is MDS, hence has the best distance possible.
Thus the weight of the codewords must increase by the amount $n-k(m+1)+1$ in each step
of the degree.
However, it is theoretically possible that certain constellations of the entries of~$G$
even allow a bigger growth rate.
\\[1ex]
\begin{proof}
We will first proof that the distance cannot be bigger than $(m+1)(n-k)+1$.
For this remember from Example~\ref{E-RScodes} that $f=\prod_{l=0}^{n-k-1}(x-\alpha^l)$
is in the code generated by~$c$. Thus $f=ac$ for some $a \in A$.
Define $\hat{g}:=f\sum_{\nu=0}^m\z^{\nu}=\sum_{\nu=0}^m\z^{\nu}\sigma^{\nu}(f)\in\Azs$.
Then $\hat{g}=ag$, hence $\hat{g}\in\lideal{g}$.
Using Theorem~\ref{T-alg}(c) we derive $\v(\hat{g})\in\cC$.
Now observe that~$f$ has weight exactly~$n-k+1$ and the same is true for
$\sigma^{\nu}(f)$ since~$\sigma$ is weight preserving. Thus
we derive at $\wt\big(\v(\hat{g})\big)=(m+1)(n-k+1)$ showing that the distance is at most
this number.

As for the rest of the theorem it suffices to prove part~(2).
Indeed, the assumption $m\leq\frac{n-k}{k}$ guarantees that $n+1-k(m+1)>0$ and thus
the lower bound in~(2) is always at least $(m+1)(n-k+1)$.
As for proving~(2), we will make use of the matrices $G_{\mu,\nu}$ from
Proposition~\ref{P-coeffcode}.
Remember from Lemma~\ref{L-coeffcodedist} that $\dist(\im G_{\mu,\nu})=n-(\nu-\mu+1)k+1$.
\\
Let $u=\sum_{j=0}^tu_j\z^j\in\F[\z]^k$ be a message with $u_0\not=0\not=u_t$ and no~$m$
consecutive zero coefficients.
Then the associated codeword $v:=uG$ has degree $t+m$ and length~$t+m+1$.
\\
In the case $t<m$ the codeword~$v$ reads as
\begin{equation}\label{e-vsmall}
  \begin{array}{rcl}
   v&=&{\DS\sum_{\nu=0}^t(u_{\nu},u_{\nu-1},\ldots,u_0)G_{0,\nu}\z^\nu
        +\sum_{\nu=t+1}^m(u_t,u_{t-1},\ldots,u_0)G_{\nu-t,\nu}\z^\nu}\\[1ex]
    & &{\DS+\sum_{\nu=m+1}^{m+t}(u_t,u_{t-1},\ldots,u_{\nu-m})G_{\nu-t,m}\z^\nu}.
  \end{array}
\end{equation}
Using Lemma~\ref{L-coeffcodedist} and the fact that $u_0\not=0\not=u_t$, we obtain for
the weight of~$v$
\begin{align*}
  \wt(v)&\geq\sum_{\nu=0}^t(n+1-k(\nu+1))+\!\!\sum_{\nu=t+1}^m\!\!(n+1-k(t+1))
            +\!\!\sum_{\nu=m+1}^{m+t}\!\!(n+1-k(m+t-\nu+1))\\
        &=(m+t+1)(n+1)-k\sum_{\nu=0}^t(\nu+1)-(m-t)k(t+1)-k\sum_{\nu=1}^{t}\nu\\
        &=(m+1)(n+1)+t(n+1-mk-k)-mk-k\\
        &=(m+1)(n-k+1)+t(n+1-k(m+1)).
\end{align*}
If $t\geq m$ one has
\begin{equation}\label{e-vbig}
 \begin{array}{rcl}
   v&=&{\DS\sum_{\nu=0}^{m-1}(u_{\nu},u_{\nu-1},\ldots,u_0)G_{0,\nu}\z^\nu
        +\sum_{\nu=m}^t(u_\nu,u_{\nu-1},\ldots,u_{\nu-m})G_{0,m}\z^\nu}\\[1ex]
    & &{\DS+\sum_{\nu=t+1}^{t+m}(u_t,u_{t-1},\ldots,u_{\nu-m})G_{\nu-t,m}\z^\nu}.
 \end{array}
\end{equation}
Using that $u_0\not=0\not=u_t$ and that no~$m$ consecutive coefficients of~$u$ are zero,
one obtains like in the previous case
\begin{align*}
 \wt(v)&\geq\sum_{\nu=0}^{m-1}\!\!(n+1-k(\nu+1))+\!\!\sum_{\nu=m}^{t}\!\!(n+1-k(m+1))+
        \!\!\sum_{\nu=t+1}^{t+m}\!\!(n+1-k(m+t-\nu+1))\\
       &=(m+t+1)(n+1)-k\sum_{\nu=0}^{m-1}(\nu+1)-(t-m+1)k(m+1)-k\sum_{\nu=1}^{m}\nu\\
       &=(m+1)(n-k+1)+t(n+1-k(m+1)).
\end{align*}
This proves the assertions.
\end{proof}
%%%%%%%%%%%%%%%%%%%%%%%%%%%%%%%%%%%%%%
\begin{rem}\label{R-biggerfield}
The results of Theorem~\ref{T-dRS} are also true if we choose the field size~$q$ such that
$n|(q-1)$ rather than $n=q-1$.
In this case there exists an element of order~$n$ in $\F$ and
this is all what is needed for the construction to work.
However, that construction does not give us a better distance and thus the constructed
codes might be farther away from the corresponding Griesmer bound for codes with parameters
$(n,k,km)$ and memory~$m$ over~$\F_q$.
\end{rem}
%%%%%%%%%%%%%%%%%%%%%%%%%%%%%%%%%%%%%%%

The following examples illustrate these results.

%%%%%%%%%%%%%%%%%%%%%%%%%%%%%%%%%%%%%%%%
\begin{exa}\label{E-RSdist}
We choose $\F=\F_8$ with primitive element~$\alpha$ satisfying $\alpha^3+\alpha+1=0$.
Thus $n=7$.
\begin{alphalist}
\item If we pick $k=2$, then the automorphism is given by $\sigma(x)=\alpha^2x$.
    The set~$S:=\{\ve{5},\ve{6}\}$, see\eqnref{e-c2}, satisfies
    $\sigma(S)=\{\ve{3},\ve{4}\},\,\sigma^2(S)=\{\ve{1},\ve{2}\},\,
    \sigma^3(S)=\{\ve{6},\ve{0}\}$.
    This shows that $b=\lfloor\frac{n}{k}\rfloor-1=2$ is the maximum value
    satisfying\eqnref{e-S}.
    We obtain
    \[
    c:=\ve{5}+\ve{6}=
    \alpha x^6+\alpha^2x^5+\alpha^2x^4+\alpha^4x^3+\alpha x^2+\alpha^4x.
    \]
    Choosing $m=2$ and applying\eqnref{e-G} we derive

    \[
    G=\begin{pmatrix}
    0&\alpha +\alpha z+\alpha z^2\\
    \alpha^4+\alpha^6 z+\alpha z^2&0\\
    \alpha+\alpha^5 z+\alpha^2 z^2&\alpha^4+\alpha z+\alpha^5 z^2\\
    \alpha^4+\alpha^3 z+\alpha^2 z^2&\alpha+z+\alpha^6 z^2\\
    \alpha^2+\alpha^3 z+\alpha^4 z^2&\alpha^4+\alpha^5 z+\alpha^6 z^2\\
    \alpha^2+\alpha^5 z+\alpha z^2&\alpha^2+\alpha^5 z+\alpha z^2\\
    \alpha+\alpha^6 z+\alpha^4 z^2&\alpha^2+z+\alpha^5 z^2
    \end{pmatrix}^{\sf T}
    \in\F[\z]^{2\times7}.
    \]
    According to Theorem~\ref{T-dRS} the code $\im G$ has distance~$18$.
    Its extended row distances satisfy $\hat{d}^r_j\geq 2j+12$ for $j\geq3$.
\item If we choose $k=3$, then the automorphism is given by 
    %%Heide: hier waren verwirrende Klammern
    $\sigma(x)=\alpha^3 x$
    %%endHeide
    and we get $S=\{\ve{4},\ve{5}\,\ve{6}\},\,\sigma(S)=\{\ve{1},\ve{2}\,\ve{3}\}$
    and $b=1$.
    In this case
    \[
    c=\ve{4}+\ve{5}+\ve{6}=\alpha^2x^6+\alpha^4x^5+\alpha^3x^4+\alpha x^3+
    \alpha^5x^2+\alpha^6x+1
    \]
    and picking $m=1$ we have
    \[
    G=\begin{pmatrix}1+z& \alpha^6+\alpha^2 z& \alpha^5+\alpha^4 z& \alpha+\alpha^3 z&
    \alpha^3+\alpha z& \alpha^4+\alpha^5 z& \alpha^2+\alpha^6 z\\
    \alpha^2+\alpha^2 z& 1+\alpha^3 z& \alpha^6+\alpha^5 z& \alpha^5+z& \alpha+\alpha^6 z&
    \alpha^3+\alpha^4 z& \alpha^4+\alpha z\\
    \alpha^4+\alpha^4 z& \alpha^2+\alpha^5 z& 1+\alpha^6 z& \alpha^6+\alpha z&
    \alpha^5+\alpha^3 z& \alpha+\alpha^2 z& \alpha^3+z
    \end{pmatrix}
    \]
    The code $\im G$ has distance~$10$ and the extended row distances satisfy
    $\hat{d}^r_j\geq 2j+6$ for $j\geq2$.
\end{alphalist}
\end{exa}
%%%%%%%%%%%%%%%%%%%%%%%%%%%%%%%%%%%%%%%%

As can be seen from the second example the rows of the generator matrices~$G$ do not
necessarily have minimal weight $(m+1)(n-k+1)$.
Indeed, since multiplication by~$x$ as well as~$\sigma$ are weight preserving maps,
each row of~$G$ has weight $(m+1)\wt(c)$ and the weight of the idempotent generator~$c$
is in general bigger than the distance~$n-k+1$ of the code~$\ideal{c}$.
This is also the reason why we have used a different element in the first paragraph
of the proof above.
Using the same idea we can actually present a generator matrix of our Reed-Solomon
convolutional codes where each row has weight $(m+1)(n-k+1)$.
Indeed, as we have seen in the first part of the proof of Theorem~\ref{T-dRS}, the
polynomial $\hat{g}=f\sum_{\nu=0}^m\z^{\nu}$ is in the left ideal $\lideal{g}$.
Since actually $f=ac$ for some unit $a\in A$ we even have $\lideal{\hat{g}}=\lideal{g}$ in
$\Azs$.
Furthermore, $\ve{i}\hat{g}=0$ for $i=0,\ldots,n-k-1$ and $\deg\ve{i}\hat{g}=m$ for
$i=n-k,\ldots,n-1$.
Therefore, just like in the proof of Theorem~\ref{T-alg}(c) the
polynomial~$\hat{g}$ is reduced in the sense of \cite[Def.~4.9(b)]{GS04}.
Using \cite[Thm.~7.8]{GS04} we obtain $\cC=\im\hat{G}$ with
\[
  \hat{G}=\begin{pmatrix}\v(\hat{g})\\ \v(x\hat{g})\\ \vdots\\ \v(x^{k-1}\hat{g}
          \end{pmatrix}
         =\sum_{\nu=0}^m\z^{\nu}\hat{G}_{\nu} \text { where }
          \hat{G}_{\nu}=\begin{pmatrix}
              \v\big(\sigma^{\nu}(f)\big)\\
              \v\big(\sigma^{\nu}(xf)\big)\\
                \vdots\\
              \v\big(\sigma^\nu(x^{k-1}f)\big)
            \end{pmatrix}\in\F^{k\times n}.
\]
Since $\wt(f)=n-k+1$ we now have that each row of~$\hat{G}$ has weight $(m+1)(n-k+1)$.
In Example~\ref{E-RSdist}(b) above we obtain
$f=\alpha^6+\alpha^5 x+\alpha^5 x^2+\alpha^2x^3+x^4$
and
\[
  \hat{G}=\begin{pmatrix}
  \alpha^6+\alpha^6z&\alpha^5+\alpha z&\alpha^5+\alpha^4z&\alpha^2+\alpha^4z&
     1+\alpha^5z&0&0\\
   0&\alpha^6+\alpha^2z&\alpha^5+\alpha^4z&\alpha^5+z&\alpha^2+z&1+\alpha z&0\\
   0&0&\alpha^6+\alpha^5z&\alpha^5+z&\alpha^5+\alpha^3z&\alpha^2+\alpha^3z&1+\alpha^4z
  \end{pmatrix}.
\]
As discussed above each row of~$\hat{G}$ has weight~$10$.

We close this section with a comparison of our codes to another construction of cyclic
convolutional codes known in the literature.
%%%%%%%%%%%%%%%%%%%%%%%%%%%%%%%%%%%%%%%
\begin{rem}\label{R-Piret88}
In~\cite[p.~445]{Pi88} Piret presents a class of cyclic convolutional codes by constructing
a suitable parity check matrix $H:=H_0+\z H_1\in\F[\z ]^{n\times(n-k)}$ where $\F=\F_{2^m}$ for
some~$m$ such that $n|(2^m-1)$ and where $k\geq\frac{n+1}{2}$.
As one can show by some straightforward computations, the resulting codes are always cyclic with respect
to the automorphism given by $\sigma(x)=x^{n-1}$ and they have dimension~$k$.
Moreover, these codes have overall constraint length~$n-k$ and unit memory, that is, all row degrees of
a minimal generator matrix are at most~$1$.
Finally, it has been shown in~\cite[p.~446]{Pi88} that the distance is $2(n-k)+1$, which is
basically due to the fact that the block code with parity check matrix $(H_0,\,H_1)$ has distance
$2(n-k)+1$.
Notice also that, since $k\geq n-k$, each minimal generator matrix of the code $\ker H$
contains $2k-n$ constant rows, therefore the code contains an $(n,2k-n)$-block code, explaining once more
that its distance cannot be bigger than $2(n-k)+1\leq n$.
These codes are best if $n-k$ is big and the optimum is reached by taking $k=\frac{n+1}{2}$ in which case
the distance is~$n$.
In contrast to that, the codes we constructed above exist only for $k\leq\frac{n}{2}$; they are best
if~$k$ is small and even optimal for $k\leq 2$. Moreover, our codes never contain constant codewords.
\end{rem}
%%%%%%%%%%%%%%%%%%%%%%%%%%%%%%%%%%%%%%%

%%%%%%%%%%%%%%%%%%%%%%%%%%%%%%%%%%%%%%%%%%
\section{A generalization to BCH codes}
\setcounter{equation}{0}
%%%%%%%%%%%%%%%%%%%%%%%%%%%%%%%%%%%%%%%%%%
In this short section we will briefly sketch how the previous ideas can be
generalized to BCH codes.
It is clear that, in principle, the computations in the proof of
Theorem~\ref{T-dRS} can be generalized to all codes of Theorem~\ref{T-alg}.
However, the resulting formulas look much more complicated.
We restrict ourselves to presenting the following case.
%%%%%%%%%%%%%%%%%%%%%%%%%%%%%%%%%%%%%%%%
\begin{prop}\label{P-isom}
Let the data be as in\eqnref{e-S} --\eqnref{e-G} and let the codes $\cC_{\mu,\nu}$
be as in Proposition~\ref{P-coeffcode}.
Assume $\dist(\cC_{\mu,\mu+\nu})=d_{\nu}$ for all $0\leq\mu\leq\mu+\nu\leq b$.
Define
\[
     D(t):=\left\{\begin{array}{ll}
                {\DS 2\sum_{\nu=0}^{t-1}d_{\nu}+(m-t+1)d_t},&\text{ if }t=0,\ldots,m,\\[1ex]
                {\DS 2\sum_{\nu=0}^{m-1}d_{\nu}+(t-m+1)d_m},&\text{ if }t\geq m+1.
           \end{array}\right.
\]
Then $\hat{d}^r_j\geq D(j-m-1)$ for all $j\geq m+1$ and
$\dist(\cC)\geq\min\{D(0),D(1),\ldots,D(m)\}$.
\end{prop}
%%%%%%%%%%%%%%%%%%%%%%%%%%%%%%%%%%%%%%%%
The assumption that all codes $\cC_{\mu,\mu+\nu}$ have the same distance independent
of~$\mu$ is satisfied whenever the automorphism~$\sigma$ is weight-preserving.
This can be seen directly from the form of $\cC_{\mu,\mu+\nu}$ given in
Proposition~\ref{P-coeffcode}.
A special case was given in Lemma~\ref{L-coeffcodedist}.
\\[1ex]
{\sc Proof:}
We argue as in the proof of Theorem~\ref{T-dRS}.
The codewords of length $m+t+1$ look again like in\eqnref{e-vsmall} and\eqnref{e-vbig}.
That gives us in the case $t<m$
\[
  \wt(v)\geq \sum_{\nu=0}^t d_{\nu}+\sum_{\nu=t+1}^m d_t+\sum_{\nu=m+1}^{m+t}d_{m+t-\nu}
        =2\sum_{\nu=0}^{t-1}d_{\nu}+(m-t+1)d_t
\]
and in the case where $t\geq m$ we obtain
\[
  \wt(v)\geq \sum_{\nu=0}^{m-1}d_{\nu}+\sum_{\nu=m}^td_m+\sum_{t+1}^{t+m}d_{m+t-\nu}
        =2\sum_{\nu=0}^{m-1}d_{\nu}+(t-m+1)d_m.
\]
This proves the first part of the proposition. The second part follows from
$$
  \dist(\cC)\geq\min\{D(t)\mid t\in\N_0\}=\min\{D(t)\mid t=0,\ldots,m\}.
  \eqno\Box
$$

In the following example we will use a BCH block code and a weight preserving
automorphism~$\sigma$.
The distances of the resulting convolutional codes will be estimated according to the
previous proposition and compared to the Griesmer bound known for the (free) distance
of convolutional codes.

%%%%%%%%%%%%%%%%%%%%%%%%%%%%%%%%%%%%%%%
\begin{exa}\label{E-isomsigma}
We choose $\F=\F_2$ and $n=31$ along with the weight preserving automorphism given by
$\sigma(x)=x^{13}$.
Then $x^{31}-1=(x-1)\prod_{i=1}^6 \pi_i$ where
\[
  \begin{array}{lll}
  \pi_1=x^5+x^2+1, & \pi_2=x^5+x^4+x^3+x^2+1, &\pi_3=x^5+x^4+x^2+x+1,\\
  \pi_4=x^5+x^3+1, & \pi_5=x^5+x^3+x^2+x+1, & \pi_6=x^5+x^4+x^3+x+1.
  \end{array}
\]
The automorphism induces the permutation $\sigma|_{E}:E\longrightarrow E$ with cycles
\[
   \big(\ve{0}\big)\big(\ve{1},\ve{2},\ve{3},\ve{4},\ve{5},\ve{6}\big).
\]
We pick the set $S:=\{\ve{1}\}$ and $b:=5$.
We will consider the codes generated by $g:=\ve{1}\sum_{i=0}^mz^i$ for all $1\leq m\leq b$.
According to Theorem~\ref{T-alg}(c) they all have dimension~$5$.
We will compute the lower bounds for the distances using Proposition~\ref{P-isom}.
Since~$\sigma$ is weight-preserving, the codes $\cC_{\mu,\mu+\nu}$ have the same distance
as $\cC_{0,\nu}=\ideal{\ve{1}+\ldots+\sigma^{\nu}(\ve{1})}$.
Moreover, according to Proposition~\ref{P-coeffcode} they are all cyclic block codes of
dimension $5(\nu+1),\ 0\leq\nu\leq b$.
We can find a lower bound of their distances by counting the number of consecutive
zeros of these codes.
In order to do so, we notice that over $\F_{32}$ with primitive element~$\alpha$ satisfying
$\alpha^5+\alpha^2+1=0$ we have
\begin{equation}\label{e-pii}
 \begin{array}{l}
  \pi_1=(x-\alpha)(x-\alpha^2)(x-\alpha^4)(x-\alpha^8)(x-\alpha^{16}),\\[.5ex]
  \pi_2=(x-\alpha^3)(x-\alpha^6)(x-\alpha^{12})(x-\alpha^{24})(x-\alpha^{17}),\\[.5ex]
  \pi_3=(x-\alpha^5)(x-\alpha^{10})(x-\alpha^{20})(x-\alpha^{9})(x-\alpha^{18}),\\[.5ex]
  \pi_4=(x-\alpha^{15})(x-\alpha^{30})(x-\alpha^{29})(x-\alpha^{27})(x-\alpha^{23}),\\[.5ex]
  \pi_5=(x-\alpha^7)(x-\alpha^{14})(x-\alpha^{28})(x-\alpha^{25})(x-\alpha^{19}),\\[.5ex]
  \pi_6=(x-\alpha^{11})(x-\alpha^{22})(x-\alpha^{13})(x-\alpha^{26})(x-\alpha^{21}).
 \end{array}
\end{equation}
It is worth mentioning that this implies
\[
   \ve{1}=\beta(x-1)\prod_{i=2}^6\pi_i
   =\gcd\big(\text{MiPo}(\alpha^i,\F_2)\,\big|\,i=17,\ldots,31\big),
\]
thus $\ideal{\ve{1}}$ is a BCH block code. Therefore we call the codes generated by~$g$
above {\em BCH convolutional codes}.
Counting successive powers of~$\alpha$ we obtain from\eqnref{e-pii} for the distances
$d_{\nu}=\dist(\cC_{0,\nu})$
\[
  d_0\geq 16,\ d_1\geq8,\ d_2\geq 8,\ d_3\geq 3,\ d_4\geq 3,\ d_5\geq 2.
\]
Using now Proposition~\ref{P-isom} we can derive lower bounds for the distances of the
codes $\cC=\im G$ where~$G$ is as in\eqnref{e-G}.
We also compare the results with the Griesmer bound known for codes with parameters
$(31,5,5m)$ and memory~$m$
\cite[Thm.~3.4]{GS03}.
\[
  \begin{array}{c|c|c}
     m & \text{lower bound for the distance} & \text{Griesmer bound}\\[.5ex]\hline\hline
    1\phantom{\Big|} & \dist(\cC)\geq\min\{32,40\}=32 & 32\\[.5ex]\hline
    2\phantom{\Big|} & \dist(\cC)\geq\min\{48,48,56\}=48 & 48\\[.5ex]\hline
    3\phantom{\Big|} & \dist(\cC)\geq\min\{64,56,64,67\}=56 & 64\\[.5ex]\hline
    4\phantom{\Big|} & \dist(\cC)\geq\min\{80,64,72,70,73\}=64 & 80\\[.5ex]\hline
    5\phantom{\Big|} & \dist(\cC)\geq\min\{96,72,80,73,76,78\}=72 & 96
  \end{array}
\]
By computing the distances $d_0$ and~$d_1$ of the two smaller block codes~$\cC_{0,0}$
and~$\cC_{0,1}$ exactly (for instance, using Maple), one obtains
$d_0=16,\,d_1=12$, showing that these two codes attain the Griesmer bound for block
codes, see \cite[Thm.~5.2.6]{vLi99}.
Computing the Gauss-Jordan form of the generator matrix of the code
$\cC_{0,2}$ one can see that $d_2=8$.
The actual value of~$d_1$
improves the lower bounds in the table above.
Indeed, for $m=3$ we obtain $\dist(\cC)\geq64$, showing that the code is optimal with
respect to its distance, and for $m=4$ we obtain $\dist(\cC)\geq78$ which is actually
pretty close to the Griesmer bound.
For $m=5$ we get $\dist(\cC)\geq 81$ which still is relatively far below the Griesmer
bound.
\\
Using the same ideas as in the proof of Theorem~\ref{T-dRS} (see also the paragraph
right after that theorem) one also obtains a lower bound for the extended row distances
of these codes. For the code with memory~$m$ this slope is given by
$d_m=\dist(\cC_{0,m})$. For instance, for 
%%Heide: Tippfehler
memory~$m=1$ 
%%endHeide
this slope is at least~$12$.
In this case we also computed the weight distribution of the code explicitly
(see also \cite[Sec.~3]{McE98a}) and obtained, using Maple,
\begin{align*}
  W(L,W)=&\frac{31L^2W^{32}}{1-6LW^{20}-15LW^{16}-10LW^2}\\
        =& 31\Big(W^{32}L^2\!+\!W^{44}(10\!+\!15W^4\!+\!6W^8)L^3\!+\!
                        W^{56}(10\!+\!15W^4\!+\!6W^8)^2L^4
           \!+\!O(L^5)\Big)
\end{align*}
showing that the least weight of atomic codewords of length~$2$ is~$32$, the least
weight of atomic codewords of length~$3$ is~$44$ etc.
Thus the slope of the extended row distances is exactly~$12$.
\end{exa}
%%%%%%%%%%%%%%%%%%%%%%%%%%%%%%%%%%%%%%%%

The numbers $D(t)$ in Proposition~\ref{P-isom} can also be generalized to the case where
$\dist(\cC_{\mu,\mu+\nu})$ does depend on~$\mu$ and~$\nu$.
However, we omit this rather technical case.

%%%%%%%%%%%%%%%%%%%%%%%%%%%%%%%%%%%%%%%%
\section{Concluding remarks}
\setcounter{equation}{0}
%%%%%%%%%%%%%%%%%%%%%%%%%%%%%%%%%%%%%%
In this paper we defined, as a special case of doubly-cyclic codes, the class of
Reed-Solomon convolutional codes, and we determined their free distance and extended
row distances.
This shows that these codes possess, at least theoretically, a good performance.
We also showed an example of how to extend these results to BCH convolutional codes.
We did not discuss the issue of decoding for these codes.
Up to now we can only come up with an iterative decoding scheme for Reed-Solomon
convolutional codes that does not outperform the algebraic decoding of the
Reed-Solomon block code $\ideal{c}$.
This certainly needs to be investigated further.

%%%%%%%%%%%%%%%%%%%%%%%%%%%%%%%%%%%%%%%%%%%%%%%%%%%
\bibliographystyle{abbrv}
\bibliography{literatureAK,literatureLZ}
%%%%%%%%%%%%%%%%%%%%%%%%%%%%%%%%%%%%%%%%%%%%%%%%%%%

\end{document}